\documentclass[10pt]{article}
\usepackage{amsmath}
\usepackage{amscd}
\usepackage{amsfonts}
\usepackage{makeidx}  % allows for index generation
\usepackage{url}
\usepackage{fancyvrb}
\usepackage{hyperref}
\usepackage{xspace}

\newtheorem{theorem}{Theorem}

\newtheorem{lemma}[theorem]{Lemma}
\newtheorem{proposition}[theorem]{Proposition}

\newtheorem{example}[theorem]{Example}

\newtheorem{remark}{Remark}

\def\ccc{\mathbb{C}}
\def\zzz{\mathbb{Z}}

\def\pf{{\bf proof}:\ }
\def\Ker{{\rm Kernel}\, }
\def\Ext{{\rm Ext}\, }
\def\Tor{{\rm Tor}\, }
\def\Im{{\rm Image}\, }
\def\Inf{{\rm Inf}\, }
\def\Res{{\rm Res}\, }

\def\qed{$\Box$}
\newcommand{\SAGE}{{\sf SAGE}\xspace}
\newcommand{\sage}{\SAGE}

\makeindex

\begin{document}

\author{David Joyner\thanks{Dept Math, US Naval Academy, Annapolis, MD, wdj@usna.edu.}}
\title{A primer on computational group homology and
  cohomology using {\tt GAP} and \sage\thanks{Dedecated to my friend and colleague 
Tony Gaglione on the occasion of his $60^{th}$ birthday.}}
%\date{1-14-2008}
\maketitle

%\tableofcontents

\vskip .4in

These are expanded lecture notes of a series of expository 
talks surveying basic aspects of group
cohomology and homology. They were written for someone who has
had a first course in graduate algebra but no background in 
cohomology. You should know the definition of a (left) 
module over a (non-commutative) ring, what $\zzz[G]$ is 
(where $G$ is a group written multiplicatively and 
$\zzz$ denotes the integers), and some ring theory and group theory.
However, an attempt has been made to (a) keep the
presentation as simple as possible, (b) either provide
an explicit reference or proof of everything.

Several computer algebra packages are used to illustrate the 
computations, though for various reasons we have focused on the 
free, open source packages, such as {\tt GAP} \cite{Gap} and 
\sage \cite{St} (which includes {\tt GAP}). 
In particular, Graham Ellis generously allowed
extensive use of his HAP \cite{Ehap} documentation (which is sometimes copied 
almost verbatim) in the presentation below.
Some interesting work not included in this (incomplete) survey
is (for example) that of Marcus Bishop \cite{Bi},
Jon Carlson \cite{C} (in MAGMA), David Green \cite{Gr}
(in C), Pierre Guillot \cite{Gu} (in GAP, C++, and \sage),  
and Marc R\"oder \cite{Ro}.

Though Graham Ellis' {\tt HAP} package (and Marc R\"oder's add-on 
{\tt HAPcryst} \cite{Ro}) can compute comhomology and homology of some
infinite groups, the computational examples given below
are for finite groups only.

\section{Introduction}

First, some words of motivation.

Let $G$ be a group and $A$ a $G$-module\footnote{We call
an abelian group $A$ (written additively) which is 
a left $\zzz[G]$-module a {\bf $G$-module}.
}.
\index{$G$-module}

Let $A^G$ denote the largest submodule of $A$ on 
which $G$ acts trivially. Let us begin by asking ourselves
the following natural question.

{\bf Question}: Suppose $A$ is a submodule of a $G$-module
$B$ and $x$ is an arbitrary $G$-fixed element of
$B/A$. Is there an element $b$ in $B$, also
fixed by $G$, which maps onto $x$ under the quotient map? 

The answer to this question can be formulated in terms of
group cohomology. (``Yes'', if $H^1(G,A)=0$.)
The details, given below, will help motivate
the introduction of group cohomology.

Let $A_G$ is the largest quotient module of $A$ on 
which $G$ acts trivially. Next, we ask ourselves
the following analogous question.

{\bf Question}: Suppose $A$ is a submodule of a $G$-module
$B$ and $b$ is an arbitrary element of
$B_G$ which maps to $0$ under the natural map
$B_G\rightarrow (B/A)_G$. Is there an element $a$ in $a_G$
which maps onto $b$ under the inclusion map? 

The answer to this question can be formulated in terms of
group homology. (``Yes'', if $H_1(G,A)=0$.)
The details, given below, will help motivate
the introduction of group homology.

Group cohomology arises as the right higher derived 
functor for $A\longmapsto A^G$.
The {\bf cohomology groups of $G$ with coefficients in $A$}
are defined by
\index{cohomology groups of $G$ with coefficients in $A$}

\[
H^n(G,A)=\Ext_{\zzz [G]}^n(\zzz ,A).
\]
(See \S \ref{sec:H^n} below for more details.)
These groups were first introduced in 1943 by S. Eilenberg and
S. MacLane \cite{EM}.
The functor $A\longmapsto A^G$ on the category of
left $G$-modules is additive and left exact.
This implies that if 

\[
0 \rightarrow 
A {\rightarrow} B {\rightarrow}  
C {\rightarrow} 0
\]
is an exact sequence of $G$-modules 
then we have a {\bf long exact sequence of cohomology} 
\index{long exact sequence of cohomology} 

\begin{equation}
\label{eqn:LESC}
\begin{array}{c}
0 \rightarrow 
A^G {\rightarrow} B^G \rightarrow C^G \rightarrow 
H^1(G,A) \rightarrow \\
H^1(G,B) \rightarrow 
H^1(G,C) \rightarrow 
H^2(G,A) \rightarrow \dots 
\end{array}
\end{equation}

Similarly, group homology arises as the left higher derived 
functor for $A\longmapsto A_G$.
The {\bf homology groups of $G$ with coefficients in $A$}
are defined by
\index{homology groups of $G$ with coefficients in $A$}

\[
H_n(G,A)=\Tor_n^{\zzz [G]} (\zzz ,A).
\]
(See \S \ref{sec:H_n} below for more details.)
The functor $A\longmapsto A_G$ on the category of
left $G$-modules is additive and right exact.
This implies that if 

\[
0 \rightarrow 
A {\rightarrow} B {\rightarrow}  
C {\rightarrow} 0
\]
is an exact sequence of $G$-modules 
then we have a {\bf long exact sequence of homology} 
\index{long exact sequence of homology} 

\begin{equation}
\label{eqn:LESH}
\begin{array}{c}
\dots  \rightarrow 
H_2(G,C) \rightarrow 
H_1(G,A) \rightarrow 
H_1(G,B) \rightarrow \\
H_1(G,C) \rightarrow A_G \rightarrow 
B_G \rightarrow C_G \rightarrow 0.
\end{array}
\end{equation}

Here we will define both
cohomology $H^n(G,A)$ and homology $H_n(G,A)$
using projective resolutions 
and the higher derived functors 
$\Ext^n$ and $\Tor_n$. We ``compute'' these
when $G$ is a finite cyclic group. We also give various 
functorial properties,
such as corestriction, inflation, restriction, and transfer.
Since some of these cohomology
groups can be computed with the help of computer algebra systems, 
we also include some discussion of how to use computers
to compute them.
We include several applications to
group theory. 

One can also define $H^1(G,A)$, $H^2(G,A)$, \dots , by
explicitly constructing cocycles and coboundaries.
Similarly, one can also define $H_1(G,A)$, $H_2(G,A)$, \dots , by
explicitly constructing cycles and boundaries.
For the proof that these constructions yield the
same groups, see Rotman \cite{R}, chapter 10.

For the general outline, we follow \S 7 in chapter 10
of \cite{R} on homology. For some details, we follow
Brown \cite{B}, Serre \cite{S} or Weiss \cite{W}. 

For a recent expository account of 
this topic, see for example Adem \cite{A}.
Another good reference is Brown \cite{B}.

\section{Differential groups}

In this section cohomology and homology are viewed in
the same framework. This ``differential groups'' 
idea was introduced by Cartan
and Eilenberg \cite{CE}, chapter IV, and developed in 
R. Godement \cite{G}, chapitre 1, \S 2. However, we shall follow
Weiss \cite{W}, chapter 1.

\subsection{Definitions}

A {\bf differential group} is a pair $(L,d)$, $L$ an abelian
group and $d:L\rightarrow L$ a homomorphism such that $d^2=0$.
We call $d$ a {\bf differential operator}. 
\index{differential operator}
\index{differential group}
The group

\[
H(L)=\Ker(d)/\Im (d)
\]
is the {\bf derived group} of $(L,d)$. If
\index{derived group}

\[
L=\oplus_{n=-\infty}^\infty L_n
\]
then we call $L$ {\bf graded}. Suppose $d$ 
(more precisely, $d|_{L_n}$) satisfies, in addition,
for some fixed $r\not= 0$, 

\[
d:L_n\rightarrow L_{n+r},\ \ \ \ n\in\zzz.
\]
We say $d$ is {\bf compatible} with the grading provided
$r=\pm 1$. In this case, we call $(L,d,r)$ a {\bf graded differential group}. As we shall see,
\index{graded differential group}
the case $r=1$ corresponds to cohomology and the
the case $r=-1$ corresponds to homology.
Indeed, if $r=1$ then we call $(L,d,r)$ a 
(differential) {\bf group of cohomology type} 
\index{differential group of cohomology type} 
and if $r=-1$ then we call $(L,d,r)$ 
a {\bf group of homology type}. 
\index{differential group of homology type} 
Note that if $L=\oplus_{n=-\infty}^\infty L_n$ is a group of 
cohomology type then $L'=\oplus_{n=-\infty}^\infty L'_n$
is a group of homology type, where $L'_n=L_{-n}$,
for all $n\in\zzz$.

\vskip .7in

{\bf For the impatient}: 
For {\it cohomology}, we shall eventually take
$L=\oplus_n {\rm Hom}_G(X_n,A)$, where the $X_n$ form a 
chain complex (with $+1$ grading) determined by a certain type of
resolution. The group $H(L)$ is an abbreviation for
$\oplus_n \Ext_{\zzz [G]}^n(\zzz,A)$.
For {\it homology}, we shall eventually take
$L=\oplus_n \zzz \otimes_{\zzz [G]}X_n$, where the $X_n$ form a 
chain complex (with $-1$ grading) determined by a certain type of
resolution. The group $H(L)$ is an abbreviation for
$\oplus_n \Tor^{\zzz [G]}_n(\zzz,A)$.

\vskip .7in

Let $(L,d)=(L,d_L)$ and $(M,d)=(M,d_M)$ be differential 
groups (to be more precise, we should use different symbols
for the differential operators of $L$ and $M$ but, for notational
simplicity, we use the same symbol and hope the context
removes any ambiguity).
A homomorphism $f:L\rightarrow M$ satisfying 
$d\circ f=f\circ d$ will be called {\bf admissible}.
\index{admissible}
For any $n\in \zzz$, we define $nf:L\rightarrow M$ 
by $(nf)(x)=n\cdot f(x)=f(x)+\dots +f(x)$ ($n$ times).
If $f$ is admissible then so is $nf$, for any $n\in \zzz$.
An admissible map $f$ gives rise to a map of 
derived groups: define the map $f_*:H(L)\rightarrow H(M)$, by
$f_*(x +dL)=f(x)+dM$, for all $x\in L$.

\subsection{Properties}
\label{sec:properties}

Let $f$ be an admissible map as above.

\begin{enumerate}
\item
The map $f_*:H(L)\rightarrow H(M)$ is a homomorphism.

\item
If $f:L\rightarrow M$ and $g:L\rightarrow M$ are admissible, then so
is $f+g$ and we have $(f+g)_*=f_*+g_*$.

\item
If $f:L\rightarrow M$ and $g:M\rightarrow N$ are admissible then
so is $g\circ f:L\rightarrow N$ and we have
$(g\circ f)_*=g_*\circ f_*$.

\item
If 

\begin{equation}
\label{eqn:LMN0}
0 \rightarrow L \stackrel{i}{\rightarrow} 
M \stackrel{j}{\rightarrow} N  \rightarrow 0
\end{equation}
is an exact sequence of differential groups with
admissible maps $i,j$ then there is a
homomorphism $d_*:H(N)\rightarrow H(L)$
for which the following triangle is exact:

{\footnotesize{
\begin{equation}
\label{eqn:LMN1}
{\footnotesize{
\begin{picture}(200.00,130.00)(-60.00,0.00)
\thicklines
\put(-30.00,50.00){$H(N)$} 
\put(5.00,65.00){\vector(1,1){50.00}} 
\put(55.00,-10.00){\vector(-1,1){50.00}} 
\put(60.00,120.00){$H(L)$} 
\put(60,-30.00){$H(M)$} 
\put(70.00,110.00){\vector(0,-1){115.00}} 
\put(80.00,60.00){$i_*$} 
\put(20.00,105.00){$d_*$} 
\put(20.00,-5.00){$j_*$} 
\end{picture}
}}
\end{equation}
}}

\vskip .7in

\noindent
This diagram\footnote{This is a special case of
Th\'eor\`eme 2.1.1 in \cite{G}.} encodes both the long exact sequence of 
cohomology (\ref{eqn:LESC}) and the long exact sequence of homology
(\ref{eqn:LESH}).

Here is the construction of $d_*$:

Recall $H(N)=\Ker(d)/\Im (d)$, so any $x\in H(N)$ is represented by
an $n\in N$ with $dn=0$. Since $j$ is surjective,
there is an $m\in M$ such that $j(m)=n$. Since $j$ is admissible
and the sequence is exact,
$j(dm)=d(j(m))=dn=0$, so $dm \in \Ker(j)=\Im (i)$.
Therefore, there is an $\ell \in L$ such that 
$dm=i(\ell)$. 
Define $d_*(x)$ to be the class
of $\ell$ in $H(L)$, i.e., $d_*(x)=\ell + dL$.

Here's the verification that $d_*$ is well-defined:

We must show that if we defined instead 
$d_*(x)=\ell' + dL$, some $\ell' \in L$, then
$\ell'-\ell\in dL$.
Pull back the above $n\in N$ 
with $dn=0$ to an $m\in M$ such that
$j(m)=n$. As above, there is an $\ell \in L$ such that 
$dm=i(\ell)$. 
Represent $x\in H(N)$ by an $n'\in N$, so $x=n'+dN$
and $dn'=0$. Pull back this $n'$ 
to an $m'\in M$ such that $j(m')=n'$. As above, there is an 
$\ell' \in L$ such that $dm'=i(\ell')$. 
We know $n'-n\in dN$, so $n'-n=dn''$, some $n''\in N$.
Let $j(m'')=n''$, some $m''\in M$, so 
$j(m'-m-dm'')=n'=n-j(dm'')=n'-n-dj(m'')=n'-n-dn''=0$.
Since the sequence $L-M-N$ is exact,
this implies there is an $\ell_0\in L$ such that 
$i(\ell_0)=m'-m-dm''$. But
$di(\ell_0)=i(d\ell_0)=dm'-dm=i(\ell')-i(\ell)=i(\ell'-\ell)$,
so $\ell'-\ell\in dL$.

\item
If $M=L\oplus N$ then $H(M)=H(L)\oplus H(N)$.

%{\footnotesize{
\pf
To avoid ambiguity, for the moment, let 
$d_X$ denote the differential operator on $X$,
where $X\in \{L,M,N\}$. 
In the notation of (\ref{eqn:LMN0}), 
$j$ is projection and $i$ is inclusion.
Since both are admissible, we know that 
$d_M|_L=d_L$ and  $d_M|_N=d_N$.
Note that $H(X)\subset X$, for any differential group $X$,
so $H(M)=H(M)\cap L\oplus H(M)\cap N\subset H(L)\oplus H(N)$.
It follows from this that that $d_*=0$. From 
the exactness of the triangle (\ref{eqn:LMN1}), 
it therefore follows that
this inclusion is an equality.
%}}

\qed

\item
Let $L$, $L'$,  $M$, $M'$,  $N$, $N'$ be differential groups.   
If

{\footnotesize{
\begin{equation}
\label{eqn:LMN}
\begin{CD}
0 @>>> L @>i>> M @>j>> N @>>> 0\\
@. @VfVV @VgVV @VhVV @. \\
0 @>>> L' @>i'>> M' @>j'>> N' @>>> 0
\end{CD}
\end{equation}
}}

\vskip .7in
\noindent
is a commutative diagram of exact sequences with
$i,i',j,j',f,g,h$ all admissible then

\[
\begin{CD}
H(L) @>i_*>> H(M) \\
@Vf_*VV @Vg_*VV  \\
H(L') @>i'_*>> H(M')
\end{CD}
\]
commutes, 

\[
\begin{CD}
H(M) @>j_*>> H(N) \\
@Vg_*VV @Vh_*VV  \\
H(M') @>i'_*>> H(N')
\end{CD}
\]
commutes, and 

\[
\begin{CD}
H(N) @>d_*>> H(L) \\
@Vh_*VV @Vf_*VV  \\
H(N') @>d_*>> H(L')
\end{CD}
\]
commutes.

%{\footnotesize{
This is a case of Theorem 1.1.3 in \cite{W} and of
Th\'eor\`eme 2.1.1 in \cite{G}.

The proofs that the first two squares commute are similar,
so we only verify one and leave the other to the reader.
By assumption, (\ref{eqn:LMN}) commutes and all the maps 
are admissible. Representing $x\in H(M)$ by
$x=m+dM$, we have

\[
\begin{split}
h_*j_*(x)&=h_*(j(m)+dN)=hj(m)+dN'=gi'(m)+dN'\\
&=
g_*(i'(m)+dM')=g_*i'_*(m+dM)=g_*i'_*(x),
\end{split}
\]
as desired.

The proof that the last square commutes is a little
different than this, so we prove this too. Represent $x\in H(N)$ by
$x=n+dN$ with $dn=0$ and recall that $d_*(x)=\ell+dL$,
where $dm=i(\ell)$,
$\ell \in L$, where $j(m)=n$, for $m\in M$.
We have

\[
f_*d_*(x)=f_*(\ell+dL)=f(\ell)+dL'.
\]
On the other hand, 

\[
d_*h_*(x)=d_*(h(n)+dN')=\ell'+dL',
\]
for some $\ell'\in L'$. Since $h(n)\in N'$, 
by the commutativity of (\ref{eqn:LMN})
and the definition of $d_*$, 
$\ell'\in L'$ is an element 
such that $i'(\ell')=gi(\ell)$. Since $i'$ is
injective, this condition on $\ell'$ determines it
uniquely mod $dL'$. By the commutativity of (\ref{eqn:LMN}), 
we may take $\ell'=f(\ell)$.
%}}

\item
Let $L$, $L'$,  $M$, $M'$,  $N$, $N'$ be differential
graded groups with grading $+1$ (i.e., of ``cohomology type'').   
Suppose that we have a commutative diagram, with all
maps admissible and all rows exact as in 
(\ref{eqn:LMN}). Then the following diagram is commutative and 
has exact rows:

{\tiny{
\[
\begin{CD}
\dots @>>> H_{n-1}(N) @>d_*>> H_n(L) @>i_*>> H_n(M) @>j_*>> H_n(N) @>d_*>> H_{n+1}(L) @>>> \dots \\
@. @Vh_*VV @Vf_*VV  @Vg_*VV  @Vh_*VV @Vf_*VV @.  \\
\dots @>>> H_{n-1}(N') @>d_*>> H_n(L') @>i'_*>> H_n(M') @>j'_*>> H_n(N') @>d_*>> H_{n+1}(L') @>>> \dots 
\end{CD}
\]
}}

This is Proposition 1.1.4 in \cite{W}. 
As pointed out there, it is an immediate consequence of the 
properties, 1-6 above.

Compare this with Proposition 10.69 in \cite{R}.

\item
Let $L$, $L'$,  $M$, $M'$,  $N$, $N'$ be differential
graded groups with grading $-1$ (i.e., of ``homology type'').   
Suppose that we have a commutative diagram, with all
maps admissible and all rows exact, as in 
(\ref{eqn:LMN}). Then the following diagram is commutative and 
has exact rows:

{\tiny{
\[
\begin{CD}
\dots @>>> H_{n+1}(N) @>d_*>> H_n(L) @>i_*>> H_n(M) @>j_*>> H_n(N) @>d_*>> H_{n-1}(L) @>>> \dots \\
@. @Vh_*VV @Vf_*VV  @Vg_*VV  @Vh_*VV @Vf_*VV @.  \\
\dots @>>> H_{n+1}(N') @>d_*>> H_n(L') @>i'_*>> H_n(M') @>j'_*>> H_n(N') @>d_*>> H_{n-1}(L') @>>> \dots 
\end{CD}
\]
}}

This is the analog of the previous property and is proven similarly.

Compare this with Proposition 10.58 in \cite{R}.

\item
Let $(L,d)$ be a differential
graded group with grading $r$.   
If $d_n=d|_{L_n}$ then $d_{n+r}\circ d_n=0$ and

\begin{equation}
\label{eqn:d_n}
\dots  \rightarrow 
L_{n-r} \stackrel{d_{n-r}}{\rightarrow} 
L_n \stackrel{d_{n}}{\rightarrow}  L_{n+r} 
\stackrel{d_{n}}{\rightarrow}  L_{n+2r} \rightarrow \dots 
\end{equation}
is exact.

\item 
If $\{L_n\ |\ n\in\zzz\}$ is a sequence of abelian groups 
with homomorphisms $d_n$ satisfying (\ref{eqn:d_n}) then
$(L,d)$ is a differential group, where $L=\oplus_n L_n$
and $d=\oplus_n d_n$.

\end{enumerate}

\subsection{Homology and cohomology}

When $r=1$, we call $L_n$ the {\bf group of $n$-cochains},
$Z_n=L_n\cap \Ker(d_n)$ the group of {\bf $n$-cocycles},
and $B_n=L_n\cap d_{n-1}(L_{n-1})$ the group of {\bf $n$-coboundaries}.
We call $H_n(L) =Z_n/B_n$ the {\bf $n^{th}$ cohomology group}.
When $r=-1$, we call $L_n$ the {\bf group of $n$-chains},
$Z_n=L_n\cap \Ker(d_n)$ the group of {\bf $n$-cycles},
and $B_n=L_n\cap d_{n+1}(L_{n+1})$ the group of {\bf $n$-boundaries}.
We call $H_n(L) =Z_n/B_n$ the {\bf $n^{th}$ homology group}.
\index{group of $n$-cochains}
\index{group of $n$-cocycles}
\index{group of $n$-cycles}
\index{group of $n$-chains}

\section{Complexes}

We introduce complexes in order to define explicit
differential groups which will then be used
to construct group (co)homology.

\subsection{Definitions}
\label{sec:complexes}

Let $R$ be a non-commutative ring, for example $R=\zzz [G]$.

We shall define a ``finite free, acyclic, augmented chain 
complex'' of left $R$-modules.

A {\bf complex} (or chain complex or $R$-complex with a negative grading) 
is a sequence of maps
\index{complex}

\begin{equation}
\label{eqn:del_n}
\dots  \rightarrow 
X_{n+1} \stackrel{\partial_{n+1}}{\rightarrow} 
X_{n} \stackrel{\partial_{n}}{\rightarrow}  X_{n-1} 
\stackrel{\partial_{n-1}}{\rightarrow}  X_{n-2} \rightarrow \dots 
\end{equation}
for which $\partial_n\partial_{n+1}=0$, for all $n$.
If each $X_n$ is a free $R$-module with a finite basis over $R$
(so is $\cong R^k$, for some $k$) then the complex is called
{\bf finite free}.
If this sequence is exact then it is called an
{\bf acyclic complex}. The complex is {\bf augmented}
if there is a surjective $R$-module homomorphism
$\epsilon : X_0\rightarrow \zzz$ and an injective
$R$-module homomorphism
$\mu : \zzz\rightarrow X_{-1}$ such that 
$\partial_0= \mu\circ \epsilon$, where (as usual) $\zzz$ is
regarded as a trivial $R$-module.
\index{acyclic complex}
\index{augmented complex}

The {\bf standard diagram} for such an $R$-complex is
\index{standard diagram}

{\footnotesize{
\[
\begin{CD}
\dots  @>>> X_2 @>\partial_2>> X_1 @>\partial_1>> X_0 @>\partial_0>> X_{-1} @>\partial_{-1}>> X_{-2} @>>> \dots \\
@. @. @. @V\epsilon VV @AA\mu A @.  \\
@. @. @. \zzz @= \zzz @.  \\
@. @. @. @VVV @AAA @.  \\
@. @. @. 0 @. 0 @.  
\end{CD}
\]
}}
Such an acyclic augmented complex can be broken up into the
{\bf positive part}

\[
\dots  \rightarrow 
X_{2} \stackrel{\partial_{2}}{\rightarrow} 
X_{1} \stackrel{\partial_{1}}{\rightarrow}  X_{0} 
\stackrel{\epsilon}{\rightarrow}  \zzz \rightarrow 0,
\]
and the {\bf negative part}

\[
0 \rightarrow 
\zzz \stackrel{\mu}{\rightarrow} 
X_{-1} \stackrel{\partial_{-1}}{\rightarrow}  X_{-2} 
\stackrel{\partial_{-2}}{\rightarrow} X_{-3} \rightarrow \dots \ .
\]
Conversely, given a positive part and a negative part,
they can be combined into a standard diagram by taking
$\partial_0=\mu\circ\epsilon$.

If $X$ is any left $R$-module, let $X^*={\rm Hom}_R(X,\zzz)$
be the dual $R$-module, where
$\zzz$ is regarded as a trivial $R$-module. 
Associated to any $f\in {\rm Hom}_R(X,Y)$ is the pull-back
$f^*\in {\rm Hom}_R(Y^*,X^*)$. (If $y^*\in Y^*$ then 
define $f^*(y^*)$ to be $y^*\circ f:X\rightarrow \zzz$.)
Since ``dualizing'' reverses the direction of the maps,
if you dualize the entire complex with a $-1$ grading, 
you will get a complex with a $+1$ grading.
This is the {\bf dual complex}.
\index{dual complex}

When $R=\zzz [G]$ then we call a
finite free, acyclic, augmented chain 
complex of left $R$-modules, a {\bf $G$-resolution}.
\index{$G$-resolution}
The maps $\partial_i:X_i\rightarrow X_{i-1}$ are sometimes
called {\bf boundary maps}.
\index{boundary maps}

\begin{remark}
{\rm
Using the command {\tt BoundaryMap}
in the {\tt GAP} {\tt CRIME} package of Marcus Bishop, one can easily compute
the boundary maps of a cohomology object associated to a $G$-module.
However, $G$ must be a $p$-group.
}
\end{remark}

\begin{example}
\label{ex:hap1}
{\rm
We use the package {\tt HAP} \cite{Ehap} to illustrate some of these
concepts more concretely.
Let $G$ be a finite group, whose elements we have ordered
in some way: $G=\{g_1,...,g_n\}$.

Since a $G$-resolution $X_*$ determines a sequence of finitely
generated free $\zzz[G]$-modules, to concretely describe $X_*$ we must 
be able to concretely describe a finite free $\zzz[G]$-module.  
In order to represent a word $w$ in a free $\zzz[G]$-module $M$ of rank $n$, 
we use a list of integer pairs $w=[ [i_1,e_1], [i_2,e_2], ..., [i_k,e_k] ]$. 
The integers $i_j$ lie in the range $\{-n,..., n\}$ and 
correspond to the free $\zzz[G]$-generators of $M$ and their 
additive inverses. The integers $e_j$ are positive (but not necessarily 
distinct) and correspond to the group element $g_{e_j}$. 

Let's begin with a {\tt HAP} computation.

\vskip .2in

\begin{Verbatim}[fontsize=\scriptsize,fontfamily=courier,fontshape=tt,frame=single,label={\tt GAP}]

gap> LoadPackage("hap");
true
gap> G:=Group([(1,2,3),(1,2)]);;
gap> R:=ResolutionFiniteGroup(G, 4);;

\end{Verbatim}

\vskip .1in
\noindent
This computes the first $5$ terms of a $G$-resolution ($G=S_3$)

\[
X_4 \stackrel{\delta_4}{\rightarrow} X_3  
 \stackrel{\delta_3}{\rightarrow} X_2
 \stackrel{\delta_2}{\rightarrow} X_1  
 \stackrel{\delta_1}{\rightarrow} X_0  
\rightarrow \zzz \rightarrow 0.
\]
The bounday maps $\delta_i$ are determined from
the {\tt boundary} component of the {\tt GAP} record {\tt R}.
This record has (among others) the following components:

\begin{itemize}
\item
{\tt R!.dimension(k)} -- the $\zzz[G]$-rank of the module $X_k$,
\item
{\tt R!.boundary(k, j)} -- the image in $X_{k-1}$ of the $j$-th free 
generator of $X_k$,
\item
{\tt R!.elts} -- the elements in $G$,
\item
{\tt R!.group} is the group in question.
\end{itemize}

Here is an illustration:

\vskip .2in

\begin{Verbatim}[fontsize=\scriptsize,fontfamily=courier,fontshape=tt,frame=single,label={\tt GAP}]

gap> R!.group;
 Group([ (1,2), (1,2,3) ])
gap> R!.elts;
 [ (), (2,3), (1,2), (1,2,3), (1,3,2), (1,3) ]
gap> R!.dimension(3);
 4
gap> R!.boundary(3,1);
 [ [ 1, 2 ], [ -1, 1 ] ]
gap> R!.boundary(3,2);
 [ [ 2, 2 ], [ -2, 4 ] ]
gap> R!.boundary(3,3);
 [ [ 3, 4 ], [ 1, 3 ], [ -3, 1 ], [ -1, 1 ] ]
gap> R!.boundary(3,4);
 [ [ 2, 5 ], [ -3, 3 ], [ 2, 4 ], [ -1, 4 ], [ 2, 1 ], [ -3, 1 ] ]

\end{Verbatim}

\vskip .1in
\noindent
In other words, $X_3$ is rank $4$ as a $G$-module, with generators
$\{f_1, f_2, f_3, f_4\}$ say, and 

\[
\delta_3(f_1) = f_1g_2 - f_1g_1,
\]
\[
\delta_3(f_2) = f_2g_2 - f_2g_4,
\]
\[
\delta_3(f_3) = f_3g_4 - f_3g_1+f_1g_3-f_1g_1,
\]
\[
\delta_3(f_4) = f_2(g_1+g_3+g_5) - f_3g_3 + f_1g_4-f_3g_1.
\]

Now, let us create another resolution and compute the 
equivariant chain map between them.
Below is the complete {\tt GAP} session:

\vskip .2in

{\footnotesize{
\begin{Verbatim}[fontsize=\scriptsize,fontfamily=courier,fontshape=tt,frame=single,label={\tt GAP}]

gap> G1:=Group([(1,2,3),(1,2)]);
Group([ (1,2,3), (1,2) ])
gap> G2:=Group([(1,2,3),(2,3)]);
Group([ (1,2,3), (2,3) ])
gap> phi:=GroupHomomorphismByImages(G1,G2,[(1,2,3),(1,2)],[(1,2,3),(2,3)]);
[ (1,2,3), (1,2) ] -> [ (1,2,3), (2,3) ]
gap> R1:=ResolutionFiniteGroup(G1, 4);
Resolution of length 4 in characteristic 0 for Group([ (1,2), (1,2,3) ]) . 

gap> R2:=ResolutionFiniteGroup(G2, 4);
Resolution of length 4 in characteristic 0 for Group([ (2,3), (1,2,3) ]) . 

gap> ZP_map:=EquivariantChainMap(R1, R2, phi);
Equivariant Chain Map between resolutions of length 4 . 

gap> map := TensorWithIntegers( ZP_map);
Chain Map between complexes of length 4 . 

gap> Hphi := Homology( map, 3);
[ f1, f2, f3 ] -> [ f2, f2*f3, f1*f2^2 ]
gap> AbelianInvariants(Image(Hphi));
[ 2, 3 ]
gap> 
gap> GroupHomology(G1,3);
[ 6 ]
gap> GroupHomology(G2,3);
[ 6 ]

\end{Verbatim}
}}
\vskip .1in
\noindent
In other words, $H(\phi)$ is an isomorphism (as it should be,
since the homology is independent of the resolution choosen).

}
\end{example}

\subsection{Constructions}

Let $R=\zzz[G]$.

\subsubsection{Bar resolution}
\label{sec:bar_res}

This section follows \S 1.3 in \cite{W}.

Define a symbol $[.]$ and call it the {\bf empty cell}.
Let $X_0=R[.]$, so $X_0$ is a finite free (left) $R$-module 
whose basis has only $1$ element.
For $n>0$, let $g_1,\dots ,g_n\in G$ and define an
{\bf $n$-cell} to be the symbol $[g_1,\dots ,g_n]$.
\index{cell}
Let

\[
X_n=\oplus_{(g_1,\dots ,g_n)\in G^n} R[g_1,\dots ,g_n],
\]
where the sum runs over all ordered $n$-tuples
in $G^n$.

Define the differential operators $d_n$ and the
augmentation $\epsilon$, as $G$-module maps, by

\[
\begin{split}
\epsilon(g[.])&=1,\ \ \ \ \ \ g\in G\\
d_1([g])&=g[.]-[.],\\
d_2([g_1,g_2])&=g_1[g_2]-[g_1g_2]+[g_1],\\
 & \vdots \\
d_n([g_1,\dots ,g_n])&=g_1[g_2,\dots ,g_n]
+\sum_{i=1}^{n-1}(-1)^i[g_1,\dots ,g_{i-1},g_ig_{i+1},g_{i+2},\dots ,g_n]\\
 & \ \ \ \ \ \ +(-1)^n[g_1,\dots ,g_{n-1}],
\end{split}
\]
for $n\geq 1$. 
Note that the condition
$\epsilon(g[.])=1$ for all $g\in G$ is equivalent to 
saying $\epsilon([.])=1$. This is because $\epsilon$ is a 
$G$-module homomorphism and $\zzz$ is a trivial $G$-module,
so $\epsilon(g[.])=g\epsilon([.])=g\cdot 1=1$, where the 
(trivial) $G$-action
on $\zzz$ is denoted by a $\cdot$.

The $X_n$ are finite free $G$-modules, with the set of all 
$n$-cells serving as a basis.

\begin{proposition} 
\label{prop:bar}
With these definitions, the sequence

\[
\dots  \rightarrow 
X_{2} \stackrel{d_{2}}{\rightarrow} 
X_{1} \stackrel{d_{1}}{\rightarrow}  X_{0} 
\stackrel{\epsilon}{\rightarrow}  \zzz \rightarrow 0,
\]
is a free $G$-resolution.
\end{proposition}

Sometimes this resolution is called the 
{\bf bar resolution}\footnote{This resolution is not
the same as the resolution computed by {\tt HAP}
in Example \ref{ex:hap1}. For details on the resolution used 
by {\tt HAP}, please see Ellis \cite{E2}.}.
\index{bar resolution}
There are two other resolutions we shall consider.
One is the closely related ``homogeneous resolution''
and the other is the ``normalized bar resolution''.

This simple-looking proposition is not so simple
to prove. First, we shall show it is a complex, i.e., 
$d^2=0$. Then, and this is 
the most non-trivial part of the proof, we show that
the sequence is exact.

First, we need some definitions and a lemma.

Let $f:L\rightarrow M$ and $g:L\rightarrow M$ be
$+1$-graded admissible maps. 
We say $f$ is {\bf homotopic} to $g$ if there is
a homomorphism $D:L\rightarrow M$,
called a {\bf homotopy}, such that 
\index{homotopy}

\begin{itemize}

\item
$D_n=D|_{L_n}:L_n\rightarrow M_{n+1}$,

\item
$f-g=Dd+dD$.
\end{itemize}
If $L=M$ and 
the identity map $1:L\rightarrow L$ is homotopic to
the zero map $0:L\rightarrow L$ then the homotopy
is called a {\bf contracting homotopy for $L$}.
\index{contracting homotopy}

\begin{lemma}
If $L$ has a contracting homotopy then $H(L)=0$.
\end{lemma}

\pf 
Represent $x\in H(L)$ by $\ell\in L$ with $d\ell =0$.
But $\ell=1(\ell)-0(\ell)=dD(\ell)+Dd(\ell)=dD(\ell)$.
Since $D:L\rightarrow L$, this shows $\ell\in dL$, 
so $x=0$ in $H(L)$.
\qed

Next, we construct a contracting homotopy for the
complex $X_*$ in Proposition \ref{prop:bar}
with differential operator $d$.
Actually, we shall {\it temporarily}
let $X_{-1}=\zzz$, $X_{-n}=0$ 
and $d_{-n}=0$ for $n>1$, 
so that that the complex is infinite in both
directions. We must define $D:X\rightarrow X$
such that 

\begin{itemize}

\item
$D_{-1}=D|_{\zzz}:\zzz\rightarrow X_0$,

\item
$D_{n}=D|_{X_n}:X_n\rightarrow X_{n+1}$,

\item
$\epsilon D_{-1}=1$ on $\zzz$,

\item
$d_1D_0+D_{-1}\epsilon =1$ on $X_0$,

\item
$d_{n+1}D_n+D_{n-1}d_n=1$ in $X_n$, for $n\geq 1$.

\end{itemize}
Define 

\[
\begin{split}
D_{-n}&=0,\ \ \ \ \ \ n>1,\\
D_{-1}(1)&=[.],\\
D_0(g[.])&=[g],\\
D_n(g[g_1,\dots ,g_n])&=[g,g_1,\dots ,g_n],\ \ \ \ \ \ n>0,
\end{split}
\]
and extend to a $\zzz$-basis linearly. 

Now we must verify the desired properties.

By definition, for $m\in \zzz$,
$\epsilon D_{-1}(m)=\epsilon (m[.])=m\epsilon ([.])=m$.
Therefore, $\epsilon D_{-1}$ is the identity map on
$\zzz$.

Similarly, 
\[
\begin{split}
(d_1 D_0+D_{-1}\epsilon )(g[.])=
d_1 ([g])+D_{-1}(1) \\
=g[.]-[.]+D_{-1}(1)=g[.]-[.]+[.]=g[.].
\end{split}
\]

For the last property, we compute

\[
\begin{split}
d_{n+1} D_n(g[g_1,\dots ,g_n])
&=d_{n+1} ([g,g_1,\dots ,g_n])\\
&=g[g_1,\dots ,g_n]-[gg_1,\dots ,g_n]\\
& \ \ \ \ \ \ 
+\sum_{i=1}^{n-1}(-1)^{i-1}[g,g_1,\dots ,
g_{i-1},g_ig_{i+1},g_{i+2},\dots ,g_n]\\
 & \ \ \ \ \ \ +(-1)^{n+1}[g,g_1,\dots ,g_{n-1}],
\end{split}
\]
and

\[
\begin{split}
D_{n-1}d_{n} (g[g_1,\dots ,g_n])\\
&=D_{n-1}(gd_{n} ([g_1,\dots ,g_n]))\\
&=D_{n-1}(gg_1[g_2,\dots ,g_n]\\
& \ \ \ \ \ \ 
+\sum_{i=1}^{n-1}(-1)^ig[g_1,\dots ,g_{i-1},g_ig_{i+1},g_{i+2},\dots ,g_n]\\
 & \ \ \ \ \ \ +(-1)^{n}g[g_1,\dots ,g_{n-1}])\\
&=[gg_1,g_2,\dots ,g_n]\\
& \ \ \ \ \ \ 
+\sum_{i=1}^{n-1}(-1)^i[g,g_1,\dots ,g_{i-1},g_ig_{i+1},g_{i+2},\dots ,g_n]\\
 & \ \ \ \ \ \ +(-1)^{n}[g,g_1,\dots ,g_{n-1}].
\end{split}
\]
All the terms but one cancels, verifying that
$d_{n+1}D_n+D_{n-1}d_n=1$ in $X_n$, for $n\geq 1$.

Now we show $d^2=0$. One verifies $d_1d_2=0$ directly
(which is left to the reader). Multiply $d_kD_{k-1}+D_{k-2}d_{k-1}=1$
on the right by $d_k$ and $d_{k+1}D_{k}+D_{k-1}d_{k}=1$
on the left by $d_k$:

\[
d_kD_{k-1}d_k + D_{k-2}d_{k-1}d_k=d_k=
d_kd_{k+1}D_{k}+d_kD_{k-1}d_{k}.
\]
Cancelling like terms, the induction hypothesis
$d_{k-1}d_k=0$ implies $d_{k}d_{k+1}=0$. This 
shows $d^2=0$ and hence that the sequence in 
Proposition \ref{prop:bar} is exact.
This completes the proof of Proposition \ref{prop:bar}.
\qed

\vskip .1in

The above complex can be ``dualized'' in the sense of \S \ref{sec:complexes}. 
This dualized complex is of the form

\[
0 \rightarrow 
\zzz \stackrel{\mu}{\rightarrow} 
X_{-1} \stackrel{d_{-1}}{\rightarrow}  X_{-2} 
\stackrel{d_{-2}}{\rightarrow} X_{-3} \rightarrow \dots \ .
\]
The {\bf standard $G$-resolution} is obtained by splicing these
together.
\index{standard $G$-resolution} 

\subsubsection{Normalized bar resolution}

Define the {\bf normalized cells} by

\[
[g_1,...,g_n]^*=
\left\{
\begin{array}{cc}
[g_1,...,g_n], &{\rm if \ all\ }g_i\not= 1,\\
0, & {\rm if \ some\ }g_i= 1.
\end{array}
\right.
\]

Let $X_0=R[.]$ and 

\[
X_n=\oplus_{(g_1,\dots ,g_n)\in G^n} R[g_1,\dots ,g_n]^*,\ \ \ \ \ 
n\geq 1,
\]
where the sum runs over all ordered $n$-tuples
in $G^n$.
Define the differential operators $d_n$ and the augmentation
map exactly as for the bar resolution.

\begin{proposition} 
\label{prop:nbar}
With these definitions, the sequence

\[
\dots  \rightarrow 
X_{2} \stackrel{d_{2}}{\rightarrow} 
X_{1} \stackrel{d_{1}}{\rightarrow}  X_{0} 
\stackrel{\epsilon}{\rightarrow}  \zzz \rightarrow 0,
\]
is a free $G$-resolution.
\end{proposition}

Sometimes this resolution is called the 
{\bf normalized bar resolution}.
\index{normalized bar resolution}

\pf 
See Theorem 10.117 in \cite{R}. \qed

\subsubsection{Homogeneous resolution}

Let $X_0=R$, so $X_0$ is a finite free (left) $R$-module 
whose basis has only $1$ element.
For $n>0$, let $X_n$ denote the $\zzz$-module generated
by all $(n+1)$-tuples $(g_0,\dots ,g_n)$. Make $X_i$ into 
a $G$-module by defining the action by
$g:X_n\rightarrow X_n$ by

\[
g:(g_0,...,g_n)\longmapsto (gg_0,\dots ,gg_n),\ \ \ \ \ g\in G.
\]

Define the differential operators $\partial_n$ and the
augmentation $\epsilon$, as $G$-module maps, by

\[
\begin{split}
\epsilon (g)&=1,\\
\partial_n(g_0,\dots ,g_n)&=\sum_{i=0}^{n-1}(-1)^i
(g_0,\dots ,g_{i-1},\hat{g}_i,g_{i+1},\dots ,g_n),
\end{split}
\]
for $n\geq 1$. 

\begin{proposition} 
\label{prop:homog}
With these definitions, the sequence

\[
\dots  \rightarrow 
X_{2} \stackrel{\partial_{2}}{\rightarrow} 
X_{1} \stackrel{\partial_{1}}{\rightarrow}  X_{0} 
\stackrel{\epsilon}{\rightarrow}  \zzz \rightarrow 0,
\]
is a $G$-resolution.
\end{proposition}

Sometimes this resolution is called the 
{\bf homogeneous resolution}.
\index{homogeneous resolution}

Of the three resolutions presented here, this 
one is the most straightforward to deal with.

\pf 
See Lemma 10.114, Proposition 10.115, and Proposition 10.116
in \cite{R}. 
 \qed

\section{Definition of $H^n(G,A)$}
\label{sec:H^n}

For convenience, we briefly recall the definition of $\Ext^n$.
Let $A$ be a left $R$-module, where $R=\zzz [G]$, and let
$(X_i)$ be a $G$-resolution of $\zzz$. We define

\[
{\rm \Ext}^n_{\zzz [G]}(\zzz,A)=
\Ker(d_{n+1}^*)/\Im (d_n^*),
\]
where 

\[
d_n^*:Hom(X_{n-1},A)\rightarrow Hom(X_n,A),
\]
is defined by sending $f:X_{n-1}\rightarrow A$ to
$fd_n:X_{n}\rightarrow A$.
It is known that this is, up to isomorphism, independent of
the resolution choosen.
Recall ${\rm \Ext}^*_{\zzz [G]}(\zzz,A)$
is the right-derived functors of the right-exact
functor $A\longmapsto A^G={\rm Hom}_G(\zzz,A)$ from the category of
$G$-modules to the category of abelian groups.
We define

\begin{equation}
\label{eqn:H^ndef}
H^n(G,A)={\rm \Ext}^n_{\zzz [G]}(\zzz,A),
\end{equation}
When we wish to emphasize the dependence on the 
resolution choosen, we write $H^n(G,A,X_*)$.

For example, let $X_*$ denote the bar 
resolution in \S \ref{sec:bar_res} above.
Call $C^n=C^n(G,A)={\rm Hom}_G(X_n,A)$ the 
{\bf group of $n$-cochains of $G$ in $A$},
\index{group of $n$-cochains of $G$ in $A$}
$Z^n=Z^n(G,A)=C^n\cap \Ker(\partial)$ the group of {\bf $n$-cocycles},
\index{group of $n$-cocycles}
and $B^n=B^n(G,A)=
\partial(C^{n-1})$ the group of {\bf $n$-coboundaries}.
\index{group of $n$-coboundaries}
We call $H^n(G,A) =Z^n/B^n$ the {\bf $n^{th}$ cohomology group of $G$
in $A$}. This is an abelian group.

We call also define the cohomology group using some other resolution,
the normalized bar resolution or the homogeneous resolution for example.
If we wish to express the dependence on the resolution $X_*$ used, we
write $H^n(G,A,X_*)$. Later we shall see that, up to isomorphism, 
this abelian group is independent of the resolution.

The group $H_2(G,\zzz)$ (which is isomorphic to the algebraic dual group of 
$H^2(G,\ccc^\times)$)
is sometimes called the {\bf Schur multiplier} of $G$. Here $\ccc$ denotes
the field of complex numbers.
\index{Schur multiplier} 

We say that the group $G$ has {\bf cohomological dimension} $n$,
\index{cohomological dimension}
written $cd(G)=n$, if
$H^{n+1}(H,A)=0$ for all $G$-modules $A$ and all subgroups
$H$ of $G$, but $H^n(H,A)\not= 0$ for some such $A$ and $H$.

\begin{remark}
\begin{itemize}

\item
If $cd(G)<\infty$ then
$G$ is torsion-free\footnote{This follows from the 
fact that if $G$ is a cyclic group then $H^n(G,\zzz)\not= 0$,
discussed below.}.

\item
If $G$ is a free abelian group of finite rank
then $cd(G)=rank(G)$. 

\item
If $cd(G)=1$ then $G$ is free.
This is a result of Stallings and Swan (see for example
\cite{R}, page 885).
\end{itemize}
\end{remark}

\subsection{Computations}

We briefly discuss computer programs which compute cohomology and
some examples of known computations.

\subsubsection{Computer computations of cohomology}

{\tt GAP} \cite{Gap} can compute some cohomology groups\footnote{See 
\S 37.22 of the {\tt GAP} manual, M. Bishop's 
package {\tt CRIME} for cohomology of $p$-groups, G. Ellis'
package {\tt HAP} for group homology and cohomology of finite
or (certain) infinite groups, and M. R\"oder's {\tt HAPCryst} package
(an add-on to the {\tt HAP} package).
\sage \cite{St} computes cohomology via it's {\tt GAP} interface.}.

All the \sage commands which compute group homology or 
cohomology require that the package {\tt HAP} be loaded. 
You can do this on the command line from the main \sage directory
by typing\footnote{This is 
the current package name - change {\tt 4.4.10\_3} to whatever 
the latest version is on
\url{http://www.sagemath.org/packages/optional/} at the time you 
read this. Also, this command assumes you are using 
\sage on a machine with an internet connection.}

\verb+sage -i gap_packages-4.4.10_3.spkg+

\begin{example}
{\rm
This example uses \sage, which wraps several of the {\tt HAP}
functions.

\vskip .2in

\begin{Verbatim}[fontsize=\scriptsize,fontfamily=courier,fontshape=tt,frame=single,label=\sage]

sage: G = AlternatingGroup(5)
sage: G.cohomology(1,7)
Trivial Abelian Group
sage: G.cohomology(2,7)
Trivial Abelian Group

\end{Verbatim}

\vskip .1in
\noindent
This implies $H^1(A_5,GF(7))=H^2(A_5,GF(7))=0$.
}
\end{example}

\subsubsection{Examples}

Some example computations of a more theoretical nature.

\begin{enumerate}

\item
$H^0(G,A)=A^G$.

This is by definition.

\item
Let $L/K$ denote a Galois extension with finite Galois
group $G$. 
We have $H^1(G,L^\times)=1$.
This is often called Hilbert's Theorem 90.

See Theorem 1.5.4 in \cite{W} or Proposition 2 in \S X.1
of \cite{S}.

\item
Let $G$ be a finite cyclic group and $A$ a trivial
torsion-free $G$-module. Then $H^1(G,A)=0$.

This is a consequence of properties given in the next section.

\item
If $G$ is a finite cyclic group of order $m$ 
and $A$ is a trivial $G$-module then

\[
H^2(G,A)=A/mA
\]

This is a consequence of properties given below.

For example, $H^2(GF(q)^\times,\ccc)=0$.

\item
If $|G|=m$, $rA=0$ and $gcd(r,m)=1$, then
$H^n(G,A)=0$, for all $n\geq 1$.

This is Corollary 3.1.7 in \cite{W}.

For example, $H^1(A_5,\zzz/7\zzz)=0$.
\end{enumerate}

\section{Definition of $H_n(G,A)$}
\label{sec:H_n}

We say $A$ is {\bf projective} if the functor
\index{projective $R$-module} 
$B\longmapsto {\rm Hom}_G(A,B)$ (from the category of $G$-modules
to the category of abelian groups) is exact.
Recall, if 

\begin{equation}
\label{eqn:P_Z}
P_{\zzz}= 
\dots  \rightarrow P_2
\stackrel{d_2}{\rightarrow} 
P_1 \stackrel{d_1}{\rightarrow}  P_0 
\stackrel{\epsilon}{\rightarrow}  \zzz \rightarrow 0
\end{equation}
is a projective resolution of $\zzz$ then

\[
\Tor_n^{\zzz [G]}(\zzz,A)
=\Ker(d_n\otimes 1_A)/\Im (d_{n+1}\otimes 1_A).
\]
It is known that this is, up to isomorphism, independent of
the resolution choosen.
Recall $\Tor_*^{\zzz [G]}(\zzz,A)$
are the right-derived functors of the right-exact
functor $A\longmapsto A_G=\zzz\otimes_{\zzz [G]}A$ 
from the category of
$G$-modules to the category of abelian groups.
We define

\begin{equation}
\label{eqn:H_ndef}
H_n(G,A)=\Tor_n^{\zzz [G]}(\zzz,A),
\end{equation}
When we wish to emphasize the dependence on the 
resolution, we write $H_n(G,A,P_\zzz)$.

\begin{remark}
{\rm
If $G$ is a $p$-group, then using the command {\tt ProjectiveResolution}
in {\tt GAP}'s {\tt CRIME} package, one can easily compute
the minimal projective resolution of a $G$-module, which can be either
trivial or given as a {\tt MeatAxe}\footnote{See for example
\url{http://www.math.rwth-aachen.de/~MTX/}.} module.
}
\end{remark}

Since we can identify the functor $A\longmapsto A_G$ with
$A\longmapsto A\otimes_{\zzz[G]} \zzz$ (where $\zzz$ is considered as a trivial 
$\zzz [G]$-module), the following is another way to 
formulate this definition.

If $\zzz$ is considered as a trivial 
$\zzz [G]$-module, then a free $\zzz [G]$-resolution of 
$\zzz$ is a sequence of $\zzz [G]$-module homomorphisms

\[
... {\rightarrow} M_n {\rightarrow} M_{n-1} 
{\rightarrow} ... {\rightarrow} M_1 {\rightarrow} M_0
\]
satisfying:
\begin{itemize}
\item
 (Freeness) Each $M_n$ is a free $\zzz[G]$-module.
\item
 (Exactness) The image of $M_{n+1} {\rightarrow} M_n$ equals 
the kernel of $M_n {\rightarrow} M_{n-1}$ for all $n>0$.
\item
 (Augmentation) The cokernel of $M_1 {\rightarrow} M_0$ is 
isomorphic to the trivial $\zzz[G]$-module $\zzz$.
\end{itemize}
The maps  $M_n {\rightarrow} M_{n-1}$ are the boundary 
homomorphisms of the resolution. 
Setting $TM_n$ equal to the abelian group $M_n/G$ obtained from 
$M_n$ by killing the $G$-action, we get an induced sequence of 
abelian group homomorphisms

\[
... {\rightarrow} TM_n {\rightarrow} TM_{n-1} 
{\rightarrow} ... {\rightarrow} TM_1 {\rightarrow} TM_0
\]
This sequence will generally not satisfy the above exactness 
condition, and one defines the integral homology of $G$ to be

\[
H_n(G,\zzz) = \Ker(TM_n  {\rightarrow} TM_{n-1}) / 
\Im (TM_{n+1}  {\rightarrow} TM_n)
\]
for all $n>0$.

\subsection{Computations}

We briefly discuss computer programs which compute homology and
some examples of known computations.

\subsubsection{Computer computations of homology}

\begin{example}
{\rm
{\tt GAP} will compute the Schur multiplier
$H_2(G,\zzz)$ using the 
\newline
{\tt AbelianInvariantsMultiplier} command.
To find 
$H_2(A_5,\zzz)$, where $A_5$ is the alternating group on $5$ letters,
type

\vskip .2in

\begin{Verbatim}[fontsize=\scriptsize,fontfamily=courier,fontshape=tt,frame=single,label={\tt GAP}]

gap> A5:=AlternatingGroup(5);
Alt( [ 1 .. 5 ] )
gap> AbelianInvariantsMultiplier(A5);
[ 2 ]

\end{Verbatim}

\vskip .1in
\noindent
So, $H_2(A_5,\ccc)\cong \zzz/2\zzz$.

Here is the same computation in \sage:
\vskip .2in

\begin{Verbatim}[fontsize=\scriptsize,fontfamily=courier,fontshape=tt,frame=single,label=\sage]

sage: G = AlternatingGroup(5)
sage: G.homology(2)
Multiplicative Abelian Group isomorphic to C2

\end{Verbatim}

}
\end{example}

\begin{example}
{\rm
The \sage command {\tt poincare\_series}
returns the Poincare series of $G \pmod p$ ($p$ must be a prime). 
In other words, if you input a (finite) permutation group $G$, 
a prime $p$, and a positive integer $n$, {\tt poincare\_series(G,p,n)} 
returns a quotient of polynomials 
$f(x)=P(x)/Q(x)$ whose coefficient of $x^k$ equals the rank of 
the vector space $H_k(G,ZZ/pZZ)$, 
for all $k$ in the range $ 1\leq k \leq n$ .

\vskip .2in

{\footnotesize{
\begin{Verbatim}[fontsize=\scriptsize,fontfamily=courier,fontshape=tt,frame=single,label=\sage]

sage: G = SymmetricGroup(5)
sage: G.poincare_series(2,10)           
 (x^2 + 1)/(x^4 - x^3 - x + 1)
sage: G = SymmetricGroup(3)
sage: G.poincare_series(2,10)     
    1/(-x + 1)

\end{Verbatim}
}}

\vskip .1in
\noindent
This last one implies

\[
\dim_{GF(2)}H_k(S_2,\zzz/2\zzz)=1, 
\]
for $1\leq k\leq 10$.
}
\end{example}

\begin{example}
\label{ex:ppart}
{\rm
Here are some more examples using \sage's interface to {\tt HAP}:

\vskip .2in

{\footnotesize{
\begin{Verbatim}[fontsize=\scriptsize,fontfamily=courier,fontshape=tt,frame=single,label=\sage]

sage: G = SymmetricGroup(5)
sage: G.homology(1)
Multiplicative Abelian Group isomorphic to C2
sage: G.homology(2)
Multiplicative Abelian Group isomorphic to C2
sage: G.homology(3)
Multiplicative Abelian Group isomorphic to C2 x C4 x C3
sage: G.homology(4)
Multiplicative Abelian Group isomorphic to C2
sage: G.homology(5)
Multiplicative Abelian Group isomorphic to C2 x C2 x C2
sage: G.homology(6)
Multiplicative Abelian Group isomorphic to C2 x C2
sage: G.homology(7)
Multiplicative Abelian Group isomorphic to C2 x C2 x C4 x C3 x C5

\end{Verbatim}
}}

\vskip .1in
\noindent
The last one means that

\[
H_7(S_5,Z) = 
(\zzz/2\zzz)^2\times (\zzz/3\zzz)\times (\zzz/4\zzz)\times (\zzz/5\zzz).
\]
%}
%\end{example}
%
%\begin{example}
%{\rm
%This example uses \sage, which wraps several of the HAP
%functions.

\vskip .2in

{\footnotesize{
\begin{Verbatim}[fontsize=\scriptsize,fontfamily=courier,fontshape=tt,frame=single,label=\sage]

sage: G = AlternatingGroup(5)
sage: G.homology(1)
Trivial Abelian Group
sage: G.homology(1,7)
Trivial Abelian Group
sage: G.homology(2,7)
Trivial Abelian Group

\end{Verbatim}
}}

\vskip .1in
\noindent
This implies $H_1(A_5,\zzz)=H_1(A_5,GF(7))=H_2(A_5,GF(7))=0$.

}
\end{example}

\subsubsection{Examples}
\label{sec:homprops}

Some example computations of a more theoretical nature.

\begin{enumerate}

\item
If $A$ is a $G$-module then $\Tor_0^{\zzz[G]}(\zzz,A)=
H_0(G,A)=A_G\cong A/DA$.

%{\footnotesize{
{\bf proof}: We need some lemmas. 

Let $\epsilon :\zzz [G]\rightarrow \zzz$
be the augmentation map. This is a ring homomorphism
(but not a $G$-module homomorphism). Let $D=\Ker(\epsilon)$
denote its kernel, the {\bf augmentation ideal}.
\index{augmentation ideal}
This is a $G$-module.

\begin{lemma}
\label{lemma:Disfree}
As an abelian group, $D$ is free abelian generated by
$G-1=\{g-1\ |\ g\in G\}$.
\end{lemma}

We write this as $D=\zzz \langle G-1\rangle$.

{\bf proof of lemma}: If $d\in D$ then $d=\sum_{g\in G}m_gg$,
where $m_g\in \zzz$ and $\sum_{g\in G}m_g=0$. Thus,
$d=\sum_{g\in G}m_g(g-1)$, so $D\subset \zzz \langle G-1\rangle$.
To show $D$ is free: If $\sum_{g\in G}m_g(g-1)=0$ 
then $\sum_{g\in G}m_g g - \sum_{g\in G}m_g=0$ in $\zzz[G]$.
But $\zzz[G]$ is a free abelian group with basis $G$, so $m_g=0$
for all $g\in G$.
\qed

\begin{lemma}
$\zzz\otimes_{\zzz [G]}A= A/DA$,
where $DA$ is generated by elements
of the form $ga-a$, $g\in G$ and $a\in A$.
\end{lemma}

Recall $A_G$ denotes the largest quotient of $A$ on which $G$ 
acts trivially\footnote{Implicit in the words ``largest quotient''
is a universal property which we leave to the reader for
formulate precisely.}.

{\bf proof of lemma}: Consider the $G$-module map,
$A\rightarrow \zzz\otimes_{\zzz[G]}A$, given by
$a\longmapsto 1\otimes a$. Since $\zzz\otimes_{\zzz[G]}A$
is a trivial $G$-module, it must factor through $A_G$. The previous 
lemma implies $A_G\cong A/DA$. (In fact, the quotient map
$q:A\rightarrow A_G$ satisfies $q(ga-a)=0$ for all $g\in G$ and 
$a\in A$, so $DA\subset \Ker(q)$. By maximality of $A_G$,
$DA=\Ker(q)$. QED) So, we have maps
$A\rightarrow A_G  \rightarrow \zzz\otimes_{\zzz[G]}A$.
By the definition of tensor products, the map
$\zzz\times A\rightarrow A_G$, $1\times a\longmapsto 1\cdot aDA$,
corresponds to a map $ \zzz\otimes_{\zzz[G]}A \rightarrow A_G $
for which the composition 
$A_G  \rightarrow \zzz\otimes_{\zzz[G]}A \rightarrow A_G $
is the identity. This forces $A_G \cong \zzz\otimes_{\zzz[G]}A$.
\qed

See also \# 11 in \S \ref{sec:basisprops}.
%}}

\item
If $G$ is a finite group then
$H_0(G,\zzz)=\zzz$.

This is a special case of the example above
(taking $A=\zzz$, as a trivial $G$-module).

\item
$H_1(G,\zzz)\cong G/[G,G]$, where $[G,G]$ is the commutator
subgroup of $G$.

%{\footnotesize{
This is Proposition 10.110 in \cite{R}, \S 10.7.

\pf 
First, we {\bf claim}: $D/D^2 \cong G/[G,G]$, where
$D$ is as in Lemma \ref{lemma:Disfree}. To prove this,
define $\theta: G\rightarrow D/D^2$ by
$g\longmapsto (g-1)+D^2$. Since
$gh-1-(g-1)-(h-1)=(g-1)(h-1)$, it follows that
$\theta(gh)=\theta(g)\theta(h)$, so $\theta$ is a homomorphism.
Since $D/D^2$ is abelian and $G/[G,G]$ is the maximal
abelian quotient of $G$, we must have $\Ker(\theta)\subset [G,G]$.
Therefore, $\theta$ factors through 
$\theta': G/[G,G]\rightarrow D/D^2$, $g[G,G]\longmapsto (g-1)+D^2$.
Now, we construct an inverse. Define
$\tau:D\rightarrow G/[G,G]$ by
$g-1 \longmapsto g[G,G]$. Since $\tau (g-1 + h-1) =
g[G,G]\cdot h[G,G]=gh [G,G]$, it is not hard to see
that this is a homomorphism. We would be essentially
done (with the construction of the inverse of $\theta'$, 
hence the proof of the 
claim) if we knew $D^2\subset \Ker(\tau)$. 
(The inverse would be the composition of
the quotient $D/D^2\rightarrow D/\Ker(\tau)$ with
the map induced from $\tau$, $D/\Ker(\tau)\rightarrow G/[G,G]$.)
This follows from the 
fact that any $x\in D^2$ can be written as 
$x=(\sum_g m_g (g-1)) (\sum_h m'_h (h-1)) =
(\sum_{g,h} m_gm'_h (g-1)(h-1))$, so
$\tau(x)=\prod_{g,h} (ghg^{-1}h^{-1})^{m_gm'_h}[G,G]=[G,G]$.
 QED (claim)

Next, we show $H_1(G,\zzz)\cong D/D^2$. From 
the short exact sequence 

\[
0\rightarrow D \rightarrow \zzz[G] 
\stackrel{\epsilon}{\rightarrow} \zzz  \rightarrow 0,
\]
we obtain the long exact sequence of homology

\begin{equation}
\label{eqn:LESHG}
\begin{array}{c}
\dots   \rightarrow 
H_1(G,D) \rightarrow 
H_1(G,\zzz[G]) \rightarrow \\
H_1(G,\zzz) \stackrel{\partial}{\rightarrow} H_0(G,D) 
\stackrel{f}{\rightarrow }
H_0(G,\zzz[G]) \stackrel{\epsilon_*}{\rightarrow} H_0(G,\zzz) \rightarrow 0.
\end{array}
\end{equation}
Since $\zzz[G]$ is a free $\zzz[G]$-module, 
$H_1(G,\zzz[G])=0$. Therefore $\partial$ is injective.
By item \# 1 above (i.e., $H_0(G,A)\cong A/DA\cong A_G$,
we have 
$H_0(G,\zzz)\cong \zzz_G=\zzz$ and
$H_0(G,\zzz[G])\cong \zzz[G]/D\cong \zzz$.
By (\ref{eqn:LESHG}), $\epsilon_*$ is surjective.
Combining the last two statements, we find 
$\zzz/\Ker(\epsilon_*)\cong \zzz$.This forces $\epsilon_*$
to be injective. This, and (\ref{eqn:LESHG}), together imply
$f$ must be $0$. Since this forces $\partial$ to be 
an isomorphism, we are done.
\qed
%}}

\item
Let $G=F/R$ be a presentation of $G$, where
$F$ is a free group and $R$ is a normal subgroup
of relations. {\bf Hopf's formula} states:
$H_2(G,\zzz)\cong (F\cap R)/[F,R]$, where $[F,R]$ is the commutator
subgroup of $G$.

See \cite{R}, \S 10.7.

The group $H_2(G,\zzz)$ is sometimes called the {\bf Schur multiplier} of $G$.
\index{Schur multiplier} 

\end{enumerate}

\section{Basic properties of $H^n(G,A)$, $H_n(G,A)$}
\label{sec:basisprops}

Let $R$ be a (possibly non-commutative) ring and $A$ be
an $R$-module. We say $A$ is {\bf injective} if the functor
\index{injective $R$-module} 
$B\longmapsto {\rm Hom}_G(B,A)$ (from the category of $G$-modules
to the category of abelian groups) is exact.
\index{injective $R$-module} 
(Recall $A$ is projective if the functor
$B\longmapsto {\rm Hom}_G(A,B)$ is exact.)
We say $A$ is {\bf co-induced} if it has the form
\index{co-induced $R$-module} 
${\rm Hom}_{\zzz}(R,B)$ for some abelian group $B$. 
We say $A$ is {\bf relatively injective} if 
\index{relatively injective $R$-module} 
it is a direct factor of a co-induced $R$-module.
We say $A$ is {\bf relatively projective} if 

\[
\begin{array}{ccc}
\pi : & \zzz [G]\otimes_\zzz A 
\rightarrow & A,\\
 & x\otimes a \longmapsto & xa,
\end{array}
\]
maps a direct factor of $\zzz [G]\otimes_\zzz A $ 
isomorphically onto $A$. These are the
$G$-modules $A$ which are isomorphic to a direct
factor of the induced module $\zzz [G]\otimes_\zzz A $.
\index{relatively projective $R$-module} 
When $G$ is finite, the notions of relatively injective
and relatively projective coincide\footnote{These notions
were introduced by Hochschild \cite{Ho}.}.

\begin{enumerate}

\item
The definition of $H^n(G,A)$ does not depend on the 
$G$-resolution $X_*$ of $\zzz$ used.
 
\item
If $A$ is an projective $\zzz [G]$-module then
$H^n(G,A)=0$, for all $n\geq 1$.

This follows immediately from the definitions.

\item
If $A$ is an injective $\zzz [G]$-module then
$H_n(G,A)=0$, for all $n\geq 1$.

See also \cite{S}, \S VII.2.

\item
If $A$ is a relatively injective $\zzz [G]$-module then
$H^n(G,A)=0$, for all $n\geq 1$.

This is Proposition 1 in \cite{S}, \S VII.2.

\item
If $A$ is a relatively projective $\zzz [G]$-module then
$H^n(G,A)=0$, for all $n\geq 1$.

This is Proposition 2 in \cite{S}, \S VII.4.

\item
If $A=A'\oplus A''$ then $H^n(G,A)=H^n(G,A')\oplus H^n(G,A'')$, 
for all $n\geq 0$.
More generally, if $I$ is any indexing family
and $A=\oplus_{i\in I} A_i$ then $H^n(G,A)=\oplus_{i\in I} H^n(G,A_i)$, 
for all $n\geq 0$.

This follows from Proposition 10.81 in \S 10.6 of
Rotman \cite{R}.

\item
If 

\[
0 \rightarrow 
A {\rightarrow} B {\rightarrow}  
C {\rightarrow} 0
\]
is an exact sequence of $G$-modules 
then we have a long exact sequence of cohomology
(\ref{eqn:LESC}).
See \cite{S}, \S VII.2, and properties of the $ext$ functor
\cite{R}, \S 10.6.

\item
$A\longmapsto H^n(G,A)$ is the higher right derived functor
associated to $A\longmapsto A^G={\rm Hom}_G(A,\zzz)$ from the category of
$G$-modules to the category of abelian groups.

This is by definition.
See \cite{S}, \S VII.2, or \cite{R}, \S 10.7.

\item
If 

\[
0 \rightarrow A {\rightarrow} B {\rightarrow}  C {\rightarrow} 0
\]
is an exact sequence of $G$-modules 
then we have a long exact sequence of homology
(\ref{eqn:LESH}).
In the case of a finite group, see \cite{S}, \S VIII.1.
In general, see \cite{S}, \S VII.4, and properties
of the $\Tor$ functor in \cite{R}, \S 10.6.

\item
$A\longmapsto H_n(G,A)$ is the higher left derived functor
associated to $A\longmapsto A_G=\zzz \otimes_{\zzz [G]}A$ 
on the category of $G$-modules.

This is by definition.
See \cite{S}, \S VII.4, or \cite{R}, \S 10.7.

\item
If $G$ is a finite cyclic group then

\[
\begin{split}
H_0(G,A) &= A_G,\\
H_{2n-1}(G,A) & = A^G/NA,\\
H_{2n}(G,A) &=\Ker(N)/DA ,
\end{split}
\]
for all $n\geq 1$. 

%{\footnotesize{
To prove this, we need a lemma.

\begin{lemma}
\label{lemma:les_cyclic}
Let $G=\langle g\rangle$ be acyclic group of order $k$. Let
$M=g-1$ and $N=1+g+g^2+...+g^{k-1}$. Then

\[
\dots  \rightarrow 
\zzz [G] \stackrel{N}{\rightarrow} 
\zzz [G] \stackrel{M}{\rightarrow}  \zzz [G] \rightarrow 
\zzz [G] \stackrel{N}{\rightarrow} 
\zzz [G] \stackrel{M}{\rightarrow}  \zzz [G] 
\stackrel{\epsilon}{\rightarrow}  \zzz \rightarrow 0,
\]
is a free $G$-resolution.

\end{lemma}

{\bf proof of lemma}:
It is clearly free. Since $MN=NM=(g-1)(1+g+g^2+...+g^{k-1})=g^k-1=0$,
it is a complex. It remains to prove exactness.
Since $\Ker(\epsilon)=D=\Im (M)$, by Lemma \ref{lemma:Disfree},
this stage is exact.

To show $\Ker(M)=\Im (N)$, let $x=\sum_{j=0}^{k-1}m_jg^j\in \Ker(M)$.
Since $(g-1)x=0$, we must have $m_0=m_1=...=m_{k-1}$. This forces
$x=m_0N\in \Im (N)$. Thus $\Ker(M)\subset \Im (N)$. Clearly
$MN=0$ implies $\Im (N)\subset \Ker(M)$, so $\Ker(M)=\Im (N)$.

To show $\Ker(N)=\Im (M)$, let $x=\sum_{j=0}^{k-1}m_jg^j\in \Ker(N)$.
Since $Nx=0$, we have $0=\epsilon(Nx)=\epsilon(N)\epsilon(x)=k\epsilon(x)$,
so $\sum_{j=0}^{k-1}m_j=0$. Observe that

\[
\begin{array}{ll}
x&=m_0\cdot 1+m_1g+m_2g^2+...+m_{k-1}g^{k-1}\\
 &=(m_0-m_0g)+(m_0+m_1)g+m_2g^2+...+m_{k-1}g^{k-1}\\
 &=(m_0-m_0g)+(m_0+m_1)g-(m_0+m_1)g^2\\
 &+(m_0+m_1+m_2)g^2-(m_0+m_1+m_2)g^3+...\\
 &+(m_0+..+m_{k-1})g^{k-1}-(m_0+..+m_{k-1})g^{k}.
\end{array}
\]
where the last two terms are actually $0$. This implies
$x=-M(m_0+(m_0+m_1)g+(m_0+m_1+m_2)g^2
+...+ (m_0+..+m_{k-1})g^{k-1}\in \Im (M)$. 
Thus $\Ker(N)\subset \Im (M)$. Clearly
$NM=0$ implies $\Im (M)\subset \Ker(N)$, so $\Ker(N)=\Im (M)$.

This proves exactness at every stage.\qed

Now we can prove the claimed property.
By property 1 in \S \ref{sec:homprops}, it suffices to assume
$n>0$. Tensor the complex in Lemma \ref{lemma:les_cyclic}
on the right with $A$:

{\footnotesize{
\[
\begin{array}{cc}
\dots  \rightarrow 
&\zzz [G]\otimes_{\zzz [G]}A \stackrel{N_*}{\rightarrow} 
\zzz [G]\otimes_{\zzz [G]}A \stackrel{M_*}{\rightarrow}  
\zzz [G]\otimes_{\zzz [G]}A \stackrel{N_*}{\rightarrow} \\
&\zzz [G]\otimes_{\zzz [G]}A \stackrel{M_*}{\rightarrow}  
\zzz [G] \otimes_{\zzz [G]}A
\stackrel{\epsilon}{\rightarrow}  \zzz\otimes{\zzz [G]}A \rightarrow 0,
\end{array}
\]
}}
where the new maps are distinguished from the old maps
by adding an asterisk. By definition, 
$\zzz [G]\otimes_{\zzz [G]}A \cong A$, and by 
property 1 in \S \ref{sec:homprops}, 
$\zzz \otimes_{\zzz [G]}A \cong A/DA$.
The above sequence becomes

\[
\dots  \rightarrow 
A \stackrel{N_*}{\rightarrow} 
A \stackrel{M_*}{\rightarrow}  
A \stackrel{N_*}{\rightarrow} 
A \stackrel{M_*}{\rightarrow}  
A
\stackrel{\epsilon}{\rightarrow}  A/DA \rightarrow 0.
\]
This implies, by definition of $\Tor$,

\[
\Tor_{2n-1}^{\zzz[G]}(\zzz,A)=\Ker(M_*)/\Im (N_*)
=A^G/NA,
\]
and
\[
\Tor_{2n}^{\zzz[G]}(\zzz,A)=\Ker(N_*)/\Im (M_*)
=A[N]/DA.
\]

See also \cite{S}, \S VIII.4.1 and the Corollary in \S VIII.4.
%}}

\item
The group $H^2(G,A)$ classifies group extensions of $A$ by $G$.

This is Theorem 5.1.2 in \cite{W}. See also \S 10.2 in 
\cite{R}.

\item
If $G$ is a finite group of order $m=|G|$ then
$mH^n(G,A)=0$, for all $n\geq 1$.

This is Proposition 10.119 in \cite{R}.

\item
If $G$ is a finite group and $A$ is a finitely-generated
$G$-module then $H^n(G,A)$ is finite, for all $n\geq 1$.

This is Proposition 3.1.9 in \cite{W} and Corollary
10.120 in \cite{R}.

\item
The group $H^1(G,A)$ constructed using resolutions
is the same as the group constructed using $1$-cocycles.
The group $H^2(G,A)$ constructed using resolutions
is the same as the group constructed using $2$-cocycles.

This is Corollary 10.118 in \cite{R}.

\item
If $G$ is a finite cyclic group then

\[
\begin{split}
H^0(G,A) &= A^G,\\
H^{2n-1}(G,A) &={\rm \Ker}\, N/DA ,\\
H^{2n}(G,A) &= A^G/NA ,
\end{split}
\]
for all $n\geq 1$.
Here $N:A\rightarrow A$ is the
norm map $Na=\sum_{g\in G}ga$ and $DA$ is the 
augmentation ideal defined above (generated by elements
of the form $ga-a$).

%{\footnotesize{
\pf
The case $n=0$: By definition, 
$H^0(G,A)=\Ext^0_{\zzz[G]}(\zzz,A)={\rm Hom}_G(\zzz,A)$.
Define $\tau : {\rm Hom}_G(\zzz,A)\rightarrow A^G$ by
sending $f\longmapsto f(1)$. It is easy to see that 
this is well-defined and, in fact, injective.
For each $a\in A^G$, define $f=f_a\in {\rm Hom}_G(\zzz,A)$ by
$f(m)=ma$. This shows $\tau$ is surjective as well,
so case $n=0$ is proven.

Case $n>0$: Applying the functor ${\rm Hom}_G(*,A)$ to the
$G$-resolution in Lemma \ref{lemma:les_cyclic} to get 

{\footnotesize{
\[
\begin{array}{cc}
\dots  \leftarrow 
&{\rm Hom}_G(\zzz [G],A) \stackrel{N_*}{\leftarrow} 
{\rm Hom}_G(\zzz [G],A) \stackrel{M_*}{\leftarrow}  {\rm Hom}_G(\zzz [G],A)
\stackrel{\epsilon_*}{\leftarrow}  {\rm Hom}_G(\zzz,A) \leftarrow 0.
\end{array}
\]
}}
It is known that ${\rm Hom}_G(\zzz[G],A)\cong A$
(see Proposition 8.85 on page 583 of \cite{R}).
It follows that

\[
\begin{array}{cc}
\dots  \leftarrow 
&A \stackrel{N_*}{\leftarrow} 
A \stackrel{M_*}{\leftarrow} A
\stackrel{\epsilon_*}{\leftarrow} A^G \leftarrow 0.
\end{array}
\]
By definition of $\Ext$, for $n>0$ we have

\[
\Ext_{\zzz[G]}^{2n}(\zzz,A)=\Ker(M_*)/\Im (N_*)=A^G/NA,
\]
and

\[
\Ext_{\zzz[G]}^{2n-1}(\zzz,A)=\Ker(N_*)/\Im (M_*)=\Ker (N)/(g-1)A,
\]
where $g$ is a generator of $G$ as 
in Lemma \ref{lemma:les_cyclic}.
\qed

See also \cite{S}, \S VIII.4.1 and the Corollary in \S VIII.4.
%}}

\item
If $G$ is a finite cyclic group 
of order $m$ and $A$ is a {\it trivial} $G$-module then

\[
\begin{split}
H^0(G,A) &= A^G,\\
H^{2n-1}(G,A) &\cong A[m],\\
H^{2n}(G,A) &\cong A/mA,
\end{split}
\]
for all $n\geq 1$.

This is a consequence of the previous property.

\end{enumerate}

\section{Functorial properties}

In this section, we investigate some of the ways in
which $H^n(G,A)$ depends on $G$.

One way to construct all these in a common framework is to 
introduce the notion of a ``homomorphism of pairs''.
Let $G,H$ be groups.
Let $A$ be a $G$-module and $B$ an $H$-module.
If $\alpha :H\rightarrow G$ is a homomorphism of groups
and $\beta:A\rightarrow B$ is a homomorphism of 
$H$-modules (using $\alpha$ to regard $B$ as an $H$-module)
then we call $(\alpha,\beta)$ a 
{\bf homomorphism of pairs}, written
\index{homomorphism of pairs}

\[
(\alpha,\beta):(G,A)\rightarrow (H,B).
\]

Let $G\subset H$ be groups and $A$ an $H$-module (so, by restriction, 
a $G$-module). We say a map

\[
f_{G,H}:H^n(G,A)\rightarrow H^n(H,A),
\]
is {\bf transitive} if $f_{G_2,G_3}f_{G_1,G_2}=f_{G_1,G_2}$,
for all subgroups $G_1\subset G_2\subset G_3$.

\index{transitive}

Let $X_*$ be a $G$-resolution and $X'_*$ a $H$-resolution,
each with a $-1$ grading.
Associated to a homomorphism of groups
$\alpha :H\rightarrow G$ is a sequence of $H$-homomorphisms

\begin{equation}
\label{eqn:hom_Hom}
A_n:X'_n\rightarrow X_n, 
\end{equation}
$n\geq 0$, such that
$d_{n+1}A_{n+1}=A_nd'_{n+1}$ and $\epsilon A_0=\epsilon'$.

\begin{theorem}
\label{thrm:hom_Hom}
\begin{enumerate}

\item
If $(\alpha,\beta):(G,A)\rightarrow (G',A')$
and $(\alpha',\beta'):(G',A')\rightarrow (G'',A'')$
are homomorphisms of pairs then so
is
$(\alpha'\circ \alpha,\beta'\circ \beta):(G,A)\rightarrow (G'',A'')$.

\item
Suppose $(\alpha,\beta):(G,A)\rightarrow (G',A')$
is homomorphism of pairs, $X_*$ is a $G$-resolution, 
and $X'_*$ is a $G'$-resolution
(each infinite in both directions, with a $-1$ grading).
Let $H^n(G,A,X_*)$ denote the derived groups associated to the
differential groups ${\rm Hom}_G(X_*,A)$ with $+1$ grading.
There is a homomorphism

\[
(\alpha,\beta)_{X_*,X'_*}:H^n(G,A,X_*)\rightarrow H^n(G',A',X'_*)
\]
satisfying the following properties.

\begin{enumerate}

\item
If $G=G'$, $A=A'$, $X=X'$, $\alpha=1$ and $\beta=1$ then 
$(1,1)_{X_*,X'_*}=1$.

\item
If $(\alpha',\beta'):(G',A')\rightarrow (G'',A'')$
is homomorphism of pairs, $X''_*$ is a $G''$-resolution
then
\[
(\alpha'\circ \alpha,\beta'\circ \beta)_{X_*,X''_*}=
(\alpha',\beta')_{X'_*,X''_*}\circ (\alpha,\beta)_{X_*,X'_*}.
\]

\item
If $(\alpha,\gamma):(G,A)\rightarrow (G',A')$
is homomorphism of pairs then
\[
(\alpha,\beta+ \gamma)_{X_*,X'_*}=
(\alpha,\beta)_{X_*,X'_*}+ (\alpha,\gamma)_{X_*,X'_*}.
\]

\end{enumerate}
\end{enumerate}
\end{theorem}

\begin{remark}
For an analogous result for homology, see \S\S III.8 in 
Brown \cite{B}.
\end{remark}

\pf We sketch the proof, following
Weiss, \cite{W}, Theorem 2.1.8, pp 52-53.

(1): This is ``obvious''.

(2): Let $(\alpha,\beta):(G,A)\rightarrow (G',A')$
be a homomorphism of pairs. Using (\ref{eqn:hom_Hom}), 
we have an associated chain map

\[
\alpha^*: {\rm Hom}_{G}(X_*,A)\rightarrow {\rm Hom}_{G'}(X'_*,A')
\]
of differential groups (Brown \S III.8 in \cite{B}).
The homomorphism of cohomology groups induced by $\alpha^*$
is denoted 

\[
\alpha^*_{n,X_*,X'_*}:H^n(G,A,X_*)\rightarrow H^n(G',A',X'_*).
\]
Properties (a)-(c) follow from \S \ref{sec:properties}
and the corresponding properties of $\alpha^*$.
\qed

As the cohomology groups are independent of the
resolution used, the map
$(\alpha,\beta)_{X_*,X'_*}:H^n(G,A,X_*)\rightarrow H^n(G',A',X'_*)$
is sometimes simply denoted by

\begin{equation}
\label{eqn:hompair}
(\alpha,\beta)_{*}:H^n(G,A)\rightarrow H^n(G',A').
\end{equation}

\subsection{Restriction}

Let $X_*=X_*(G)$ denote the bar resolution.

If $H$ is a subgroup of $G$ then
the cycles on $G$, $C^n(G,A)={\rm Hom}_G(X_n(G),A)$,
can be restricted to $H$:
$C^n(H,A)={\rm Hom}_H(X_n(H),A)$. The restriction map
$C^n(G,A)\rightarrow C^n(H,A)$ leads to a
map of cohomology classes:
\[
\Res:H^n(G,A)\rightarrow H^n(H,A).
\]

In this case, the homomorphism of pairs is 
given by the inclusion map $\alpha:H\rightarrow G$ 
and the identity map $\beta:A\rightarrow A$.
The map $\Res$ is the induced map defined by 
(\ref{eqn:hompair}). By the properties of 
this induced map, we see that 
$\Res_{H,G}$ is transitive: if
$G\subset G'\subset G''$ then\footnote{There
is an analog of the restriction for
homology which also satisfies this transitive
property (Proposition 9.5 in Brown \cite{B}).}
\[
\Res_{G',G}\circ \Res_{G'',G'}=\Res_{G'',G}.
\]

A particularly nice feature of the restriction map is the 
following fact.

\begin{theorem}
If $G$ is a finite group and $G_p$ is a $p$-Sylow subgroup 
and if $H^n(G,A)_p$ is the $p$-primary component of $H^n(G,A)$
then 

(a) there is a canonical isomorphism
$H^n(G,A)\cong \oplus_p H^n(G,A)_p$, and

(b) $Res:H^n(G,A)\rightarrow H^n(G_p,A)$ 
restricted to $H^n(G,A)_p$ (identified with a subgroup
of $H^n(G,A)$ via (a)) is injective.

\end{theorem}

\pf
See Weiss, \cite{W}, Theorem 3.1.15.
\qed

\begin{example}
\label{ex:sylow}
{\rm
Homology is a functor. That is, for any $n > 0$ and group 
homomorphism $f : G \rightarrow G'$ there is an induced homomorphism  
$H_n(f) : H_n(G,\zzz) \rightarrow H_n(G',\zzz)$ satisfying

\begin{itemize}

\item
$H_n(gf) = H_n(g)H_n(f)$ for group homomorphisms $f : G \rightarrow G'$
$g : G' \rightarrow G''$,
\item
$H_n(f)$ is the identity homomorphism if $f$ is the identity.
\end{itemize}

The following commands compute $H_3(f) : H_3(P,\zzz) \rightarrow
H_3(S_5,\zzz)$ for the inclusion $f : P \hookrightarrow S_5$ 
into the symmetric group $S_5$ of its Sylow $2$-subgroup. They 
also show that the image of the induced homomorphism $H_3(f)$ 
is precisely the Sylow $2$-subgroup of $H_3(S_5,\zzz)$.

\vskip .2in

{\footnotesize{
\begin{Verbatim}[fontsize=\scriptsize,fontfamily=courier,fontshape=tt,frame=single,label={\tt GAP}]

gap>  S_5:=SymmetricGroup(5);;  P:=SylowSubgroup(S_5,2);;
gap>  f:=GroupHomomorphismByFunction(P,S_5, x->x);;
gap>  R:=ResolutionFiniteGroup(P,4);;
gap>  S:=ResolutionFiniteGroup(S_5,4);;
gap>  ZP_map:=EquivariantChainMap(R,S,f);;
gap>  map:=TensorWithIntegers(ZP_map);;
gap>  Hf:=Homology(map,3);;
gap>  AbelianInvariants(Image(Hf));
[2,4]
gap>  GroupHomology(S_5,3);
[2,4,3]

\end{Verbatim}
}}

}
\end{example}

If $H$ is a subgroup of finite index in $G$ then 
there is an analogous restriction map
in group homology (see for example Brown \cite{B}, \S III.9).

\subsection{Inflation}

Let $X_*$ denote the bar resolution of $G$. Recall
\[
X_n=\oplus_{(g_1,\dots ,g_n)\in G^n} R[g_1,\dots ,g_n],
\]
where the sum runs over all ordered $n$-tuples
in $G^n$. If $H$ is a subgroup of $G$, let $X^H_*$ denote the complex
defined by
\[
X^H_n=\oplus_{(g_1H,\dots ,g_nH)\in (G/H)^n} R[g_1H,\dots ,g_nH].
\]
This is a resolution, and we have a chain map
defined on $n$-cells by 
$[g_1,\dots ,g_n]\longmapsto [g_1H,\dots ,g_nH]$.

Suppose that $H$ is a normal subgroup of $G$
and $A$ is a $G$-module. We may view $A^H$ as a
$G/H$-module. 
In this case, the homomorphism of pairs is 
given by the quotient map $\alpha:G\rightarrow G/H$ 
and the inclusion map $\beta:A^H\rightarrow A$.
The {\bf inflation} map $\Inf$ is the induced map defined by 
(\ref{eqn:hompair}), denoted
\index{inflation} 

\[
\Inf:H^n(G/H,A^H)\rightarrow H^n(G,A).
\]

%This map is also transitive.

The {\bf inflation-restriction sequence in dimension $n$} is
\index{inflation-restriction sequence in dimension $n$}

\[
0\rightarrow H^n(G/H,A^H)\stackrel{\Inf}{\rightarrow} 
H^n(G,A) \stackrel{\Res}{\rightarrow} 
H^n(H,A).
\]
For a proof, see Weiss, \cite{W}, \S 3.4.

There an analog of this inflation-restriction 
sequence for homology.

We omit any discussion of transfer and Shapiro's lemma, due
to space limitations.

{\it Acknowledgements}:
I thank G. Ellis, M. Mazur and J. Feldvoss,
P. Guillot for correspondence
which improved the content of these notes.
%I thank Pierre Guillot for several references given below.

%\printindex

\end{document}